\documentclass[12pt,reqno]{amsart}

\usepackage{amssymb}

\usepackage{mathrsfs}

\DeclareMathAlphabet{\mathpzc}{OT1}{pzc}{m}{it}

\usepackage[all]{xy}

\entrymodifiers={+!!<0pt,\fontdimen22\textfont2>}

\def\BA{{\mathbb A}}
\def\BB{{\mathbb B}}
\def\BC{{\mathbb C}}

\def\BF{{\mathbb F}}
\def\BG{{\mathbb G}}

\def\BI{{\mathbb I}}
\def\BJ{{\mathbb J}}

\def\BQ{{\mathbb Q}}

\def\BT{{\mathbb T}}

\def\BZ{{\mathbb Z}}

\def\fa{{\mathfrak a}}
\def\fb{{\mathfrak b}}

\def\ff{{\mathfrak f}}
\def\frakg{{\mathfrak g}}

\def\fm{{\mathfrak m}}

\def\Add{\operatorname{Add}}

\def\dirlim{\operatorname{lim}_{_{_{_{\hspace{-17pt} 
                                      \textstyle \longrightarrow}}}}}

\def\Ext{\operatorname{Ext}}

\def\Hom{\operatorname{Hom}}

\def\Ker{\operatorname{Ker}}

\def\Mod{\operatorname{Mod}}
\def\opp{\operatorname{o}}

\def\prod{\operatorname{prod}}

\def\Tor{\operatorname{Tor}}

\numberwithin{equation}{part}

\newtheorem{Lemma}{Lemma}[section]
\newtheorem{Theorem}[Lemma]{Theorem}
\newtheorem{Proposition}[Lemma]{Proposition}

\theoremstyle{definition}
\newtheorem{Definition}[Lemma]{Definition}

\newtheorem{Remark}[Lemma]{Remark}

\newtheorem{Example}[Lemma]{Example}

\begin{document}

\setlength{\parindent}{0pt}
\setlength{\parskip}{6pt}

\title[Covers and purity]
{Covers, preenvelopes, and purity}

\author{Henrik Holm \ \ }
\address{Department of Mathematical Sciences, University of Aarhus, Ny
Mun\-ke\-ga\-de, Building 1530, 8000 Aarhus C, Denmark}
\email{holm@imf.au.dk}
\urladdr{http://home.imf.au.dk/holm}

\author{\ \ Peter J\o rgensen}
\address{School of Mathematics and Statistics, Newcastle University,
Newcastle upon Tyne NE1 7RU, United Kingdom}
\email{peter.jorgensen@ncl.ac.uk}
\urladdr{http://www.staff.ncl.ac.uk/peter.jorgensen}




\begin{abstract} 
  
We show that if a class of modules is closed under pure quotients,
then it is precovering if and only if it is covering, and this
happens if and only if it is closed under direct sums.  This is
inspired by a dual result by Rada and Saor\'{\i}n.

We also show that if a class of modules contains the ground ring and
is closed under extensions, direct sums, pure submodules, and pure
quotients, then it forms the first half of a so-called perfect
cotorsion pair as introduced by Salce; this is stronger than being
covering.

Some applications are given to concrete classes of modules such as
kernels of homological functors and torsion free modules in a torsion
pair. 

\end{abstract}

\maketitle

\setcounter{section}{-1}
\section{Introduction}
\label{sec:introduction}

\noindent
{\bf Covers. }
The main topic of this paper is the notion of covering classes.  To
explain what that means, observe that the classical
ho\-mo\-lo\-gi\-cal algebra of a ring can be phrased in terms of the
class of projective modules.  This class permits the construction of
projective resolutions which again enable the computation of derived
functors.

In relative homological algebra, the class of projective modules is
replaced by another, suitably chosen class of modules.  This replaces
projective resolutions by resolutions in terms of modules in the
chosen class, and derived functors by relative derived functors.  A
classical example of this is pure homological algebra where the
projective modules are replaced by the so-called pure projective
modules; these are the direct summands in direct sums of finitely
presented mo\-du\-les.  Pure homological algebra is a useful tool with
a number of applications; see for instance \cite{JL}.

Some conditions have to be imposed on a class if it is to be a
suitable replacement for the projectives, and this leads to the notion
of pre\-co\-ve\-ring classes.  These can be used instead of the projective
modules for doing homological algebra.  A class $\BF$ of modules over
a ring is precovering (or, as it is also called, contravariantly
finite) if each module $M$ has an $\BF$-precover, that is, a
homomorphism $F \rightarrow M$ with $F$ in $\BF$, which has the
property that each homomorphism $F^{\prime} \rightarrow M$ with
$F^{\prime}$ in $\BF$ has a factorization
\[
  \xymatrix{
    & F \ar[d] \\
    F^{\prime} \ar[r] \ar[ru] & M \lefteqn{.} \\
  }
\]
A precovering class enables the construction of well behaved
resolutions: Pick a precover $F_0 \rightarrow M$, let $K_0$ be the
kernel, pick a precover $F_1 \rightarrow K_0$, let $K_1$ be the
kernel, and so on.  This gives a complex
\[
  \cdots \rightarrow F_1 \rightarrow F_0 \rightarrow M \rightarrow 0
\] 
which is called a proper $\BF$-resolution of $M$.  It has the property
that it becomes exact when the functor $\Hom(F,-)$ is applied to it
for each $F$ in $\BF$.  This implies that it is unique up to homotopy,
and hence well suited for homological tasks such as the computation of
relative derived functors.

Covering classes arise as a sharpening of the notion of precovering
classes.  A class $\BF$ is covering if each module $M$ has an
$\BF$-cover $F \rightarrow M$, that is, an $\BF$-precover with the
additional property that if $F \rightarrow F$ is an endomorphism for
which $F \rightarrow F \rightarrow M$ equals $F \rightarrow M$, then
$F \rightarrow F$ is in fact an automorphism.

The notions of precovering and covering classes can be dualized, and
this results in the notions of preenveloping and enveloping classes.

Considerable energy has gone into proving that concrete classes are
(pre)covering or (pre)enveloping under suitable conditions on the
gro\-und ring.  Examples include the classes of modules which are
projective, flat, injective, Gorenstein projective, Gorenstein flat,
Gorenstein injective, pure projective, pure injective, of projective
dimension $\leq n$, torsion free, and cotorsion.  A number of these
results can be found in Enochs and Jenda's pivotal book \cite{EJ}, but
see also \cite{Adams}, \cite{ATT}, \cite{BBE}, \cite{BT}, \cite{BT2},
\cite{EklofTrlifaj}, \cite{ElBashir}, \cite{Enochs}, \cite{ELR},
\cite{KS}, \cite{Pinzon}, \cite{RadaSaorin}, and \cite{Teply2}.

This paper shows that classes possessing some simple properties from
pure homological algebra are covering, as we shall now describe.

\noindent
{\bf Purity. }
Consider a short exact sequence
\[
  0 \rightarrow M^{\prime} \rightarrow M \rightarrow M^{\prime\prime} \rightarrow 0
\]
where $M^{\prime}$ is a submodule of $M$ and $M^{\prime\prime}$ is the
corresponding quotient module.  The sequence is called {\em pure
exact} if it stays exact when tensored with any module, and then
$M^{\prime}$ is called a pure submodule and $M^{\prime\prime}$ is
called a pure quotient module of $M$.

Recall the clever result \cite[cor.\ 3.5(c)]{RadaSaorin} by Rada and
Saor\'{\i}n, that if a class $\BG$ of modules over a ring is closed
under pure submodules, then $\BG$ is preenveloping if and only if it
is closed under direct products.

Our first main result (Theorem \ref{thm:covering}) is the dual of
this.  In fact, we prove more than the dual, namely, if a class $\BF$
is closed under pure quotient modules, then $\BF$ is precovering if
and only if it is covering, and this happens if and only if $\BF$ is
closed under direct sums.  The proof uses different methods from those
of Rada and Saor\'{\i}n which do not dualize.

We go on to show that if $\BF$ contains the ground ring and is closed
under extensions, direct sums, pure submodules, and pure quotients,
then $\BF$ is the first half of a so-called perfect cotorsion pair
(Theorem \ref{thm:cotorsion}); this is a stronger property than being
covering.  Cotorsion pairs go back to Salce \cite{Salce}, and have
gained popularity as a framework for relative homological algebra.
The formal definition is stated in Definition
\ref{def:cotorsion_pairs}; the book \cite{EJ} is a useful reference,
but see also \cite{ATT}, \cite{BBE}, \cite{EklofTrlifaj}, \cite{ELR}
and \cite{KS}.

\noindent
{\bf Applications. }  As an application of these results, we
investigate classes of the form $\Ker \Ext^1(\BA,-)$, $\Ker
\Tor_1(\BB,-)$, and $\Ker \Ext^1(-,\BC)$, where $\BA$ is a class of
finitely presented modules, $\BB$ is any class of modules, and $\BC$
is a class of pure injective modules.  The notation is
straightforward; for instance,
\[
  \Ker \Tor_1(\BB,-)
  = \{\, M \,|\, \Tor_1(B,M) = 0 \; \mbox{for each $B$ in $\BB$} \,\}.
\]
Such classes had been studied previously; for instance, it was proved
by Eklof and Trlifaj in \cite[cor.\ 10 and thm.\ 12(i)]{EklofTrlifaj} that
$\Ker \Tor_1(\BB,-)$ and $\Ker \Ext^1(-,\BC)$ are both covering, and when
the ground ring is left-coherent, El Bashir's result \cite[thm.\ 
3.3]{ElBashir} implies that $\Ker \Ext^1(\BA,-)$ is also covering if
one is willing to assume Vopenka's Principle on high cardinal
numbers.

However, among other things, we prove (Theorem \ref{thm:ABC}) that if
the ground ring is left-coherent, then
\[
  \Ker \Ext^1(\BA,-) \textstyle{\bigcap} \Ker \Ext^1(-,\BC)
\]
is covering.  We also give some concrete examples of classes of this
form (see Example \ref{exa:ABC}), including the class of left-modules
of flat dimension $\leq m$ and injective dimension $\leq n$ over a
left-noetherian ring, and the class of fp-injective left-modules.

The fp-injective left-modules had already been proved to be
preenveloping by Adams \cite{Adams}, and over a left-coherent ring,
they had been proved to be covering by Pinzon \cite{Pinzon}, but we
recover their results with new proofs.  We also use our theory to give
new proofs of the following: In a so-called hereditary torsion pair of
finite type, the torsion free modules form a covering class (Theorem
\ref{thm:torsion_free_covering}), and if, moreover, the ground ring is
itself torsion free, then the torsion free modules form the first half
of a cotorsion pair (Theorem \ref{thm:torsion_and_cotorsion}).  These
results were first proved by different methods by Bican and
Torrecillas in \cite[cor.\ 4.1]{BT2}, and Angeleri-H\"{u}gel, Tonolo,
and Trlifaj in \cite[exa.\ 2.7]{ATT}.

\noindent
{\bf Notation. }
Our notation is standard and should not require explanations, but we
do wish to introduce the following blanket items.

Throughout the paper, $R$ is a ring and the word {\em class} means
{\em class of $R$-left-modules closed under isomorphisms}.

The cardinality of a module $M$ is denoted by $|M|$.

\section{Cardinality and co-cardinality classes}
\label{sec:Enochs}

In this preliminary section, we introduce the notions of cardinality
and co-cardinality classes.  They are inspired by \cite[Props.\ 5.2.2
and 6.2.1]{EJ}, and can be used to prove that classes of modules
are pre\-co\-ve\-ring and preenveloping.

\begin{Definition}
\label{def:Enochs}
\begin{enumerate}

  \item  A class $\BF$ is called a {\em cardinality class} when, for each
         $R$-left-module $M$, there is a cardinal number $\ff$ such
         that each homomorphism
\[
  F \rightarrow M
\]
         with $F$ in $\BF$ can be factored as
\[
  F \rightarrow F^{\prime} \rightarrow M
\]
         with $F^{\prime}$ in $\BF$ satisfying $|F^{\prime}| \leq
         \ff$. 

\smallskip

  \item  A class $\BG$ is called a {\em co-cardinality class} when, for
         each $R$-left-module $N$, there is a cardinal number $\frakg$
         such that each homomorphism
\[
  N \rightarrow G
\]
         with $G$ in $\BG$ can be factored as
\[
  N \rightarrow G^{\prime} \rightarrow G
\]
         with $G^{\prime}$ in $\BG$ satisfying $|G^{\prime}| \leq
         \frakg$.
\end{enumerate}
\end{Definition}

The next proposition is very close to being in \cite{EJ}, but we think
it is worth to state it explicitly.

\begin{Proposition}
\label{pro:precovering_preenveloping}
\begin{enumerate}

  \item Let $\BF$ be a class which is closed under direct summands.
        Then $\BF$ is precovering if and only if it is a cardinality class
        which is closed under set indexed direct sums.

\smallskip

  \item Let $\BG$ be a class which is closed under direct summands.
        Then $\BG$ is preenveloping if and only if it is a
        co-car\-di\-na\-li\-ty class which is closed under set indexed
        direct products.

\end{enumerate}
\end{Proposition}

\begin{proof}
It is enough to prove (i) since (ii) is dual.

If $\BF$ is a cardinality class which is closed under set indexed
direct sums, then it is precovering by \cite[prop.\ 5.2.2]{EJ}.
Conversely, if $\BF$ is pre\-co\-ve\-ring, then it is a cardinality
class, also by \cite[prop.\ 5.2.2]{EJ}.

Finally, let $\BF$ be precovering and let $\{ F_i \}$ be a set indexed
system from $\BF$.  Pick an $\BF$-precover $F \rightarrow \bigoplus
F_i$.  Each $F_j$ has an inclusion into $\bigoplus F_i$, and since $F_j$
is in $\BF$, each inclusion lifts through $F \rightarrow \bigoplus F_i$.
Taken together, this gives a splitting of $F \rightarrow \bigoplus F_i$,
so $\bigoplus F_i$ is a direct summand of $F$.  But $\BF$ is closed
under direct summands, so $\bigoplus F_i$ is in $\BF$.
\end{proof}

\begin{Example}
If $\BB$ is a set of $R$-left-modules, then it is easy to verify that
$\Add \BB$, the class of modules which are isomorphic to a direct
summand of a set indexed direct sum of modules from $\BB$, is
a co-cardinality class.

Hence, if $\Add \BB$ is closed under set indexed direct products, then
it is preenveloping by Proposition
\ref{pro:precovering_preenveloping}(ii).  Note that in some cases of
interest, $\Add \BB$ is indeed closed under set indexed direct
products, for example if $\BB = \{ B \}$ for a finitely generated
module $B$ over an Artin algebra, cf.\ \cite[lem.\ 1.2]{KS}.
\end{Example}

\section{Purity}
\label{sec:purity}

This section shows that if a class $\BF$ is closed under pure quotient
mo\-du\-les, then $\BF$ is precovering if and only if it is covering,
and this happens if and only if $\BF$ is closed under direct sums
(Theorem \ref{thm:covering}).

The following lemma is a standard application of Zorn's lemma.

\begin{Lemma}
\label{lem:maximal}
Given an inclusion of modules $K \subseteq F$, there exists
a $K^{\prime}$ which is maximal with the properties that $K^{\prime}
\subseteq K \subseteq F$ and that $K^{\prime}$ is a pure submodule of
$F$.  
\end{Lemma}

The next lemma is a special case of \cite[thm.\ 5]{BBE}.

\begin{Lemma}
\label{lem:BBE}
For each cardinal number $\fm$ there exists a cardinal number $\ff$,
depending only on $\fm$ and the ground ring $R$, such that if an
inclusion of $R$-left-modules $K \subseteq F$ has
\[
  |F/K| \leq \fm \;\; \mbox{and} \;\; |F| \geq \ff
\]
then there exists
\[
  0 \not= K^{\prime\prime} \subseteq K \subseteq F
\]
such that $K^{\prime\prime}$ is a pure submodule of $F$.
\end{Lemma}

The proof of the following proposition is inspired by the proof of
\cite[thm.\ 6]{BBE}.

\begin{Proposition}
\label{pro:Enochs}
If a class $\BF$ is closed under pure quotient modules, then it is
a cardinality class (cf.\ Definition \ref{def:Enochs}(i)).
\end{Proposition}

\begin{proof}
Let $M$ be an $R$-left-module of cardinality $\fm$ and let $\ff$ be
the cardinal number from Lemma \ref{lem:BBE}.  Let 
\[
  F \rightarrow M
\]
be a homomorphism with $F$ in $\BF$.  We will construct a
factorization as required by Definition \ref{def:Enochs}(i).

If $|F| \leq \ff$, then consider the factorization of $F
\rightarrow M$ as 
\[
  F \rightarrow F \rightarrow M,
\]
where the first arrow is the identity.  This meets the requirements of
Definition \ref{def:Enochs}(i).

If $|F| > \ff$, then let 
\[
  K = \Ker(F \rightarrow M)
\]
and use Lemma \ref{lem:maximal} to find $K^{\prime}$ maximal with the
properties that
\[
  K^{\prime} \subseteq K \subseteq F
\]
and that $K^{\prime}$ is a pure submodule of $F$.  Then $F \rightarrow
M$ has the fac\-to\-ri\-za\-ti\-on
\[
  F \rightarrow F/K^{\prime} \rightarrow M,
\]
and we will show that this meets the requirements of Definition
\ref{def:Enochs}(i).  

First, $\BF$ is closed under pure quotients, so $F/K^{\prime}$ is in
$\BF$.

Secondly, we must show $|F/K^{\prime}| \leq \ff$.  Assume to the
contrary that
\[
  |F/K^{\prime}| > \ff.
\]
Consider the inclusion
\[
  K/K^{\prime} \subseteq F/K^{\prime}.
\]
Since $F/K$ is isomorphic to a submodule of $M$, we have
\[
  \bigg| \frac{F/K^{\prime}}{K/K^{\prime}} \bigg| = |F/K| \leq |M| = \fm,
\]
and hence Lemma \ref{lem:BBE} says that there exists
\[
  0 \not= K^{\prime\prime}/K^{\prime} \subseteq K/K^{\prime} \subseteq F/K^{\prime}
\]
such that $K^{\prime\prime}/K^{\prime}$ is a pure submodule of
$F/K^{\prime}$.  We now have
\[
  K^{\prime} \subsetneqq K^{\prime\prime} \subseteq K \subseteq F, 
\]
and we claim that $K^{\prime\prime}$ is in fact a pure submodule of
$F$, contradicting the maximality of $K^{\prime}$.

For this, consider the commutative diagram
\[
  \xymatrix{
    0 \ar[r] & K^{\prime} \ar[r] \ar@{=}[d] 
             & K^{\prime\prime} \ar[r] \ar@{^{(}->}[d] 
             & K^{\prime\prime}/K^{\prime} \ar[r] \ar@{^{(}->}[d]
             & 0 \\
    0 \ar[r] & K^{\prime} \ar[r]
             & F \ar[r]
             & F/K^{\prime} \ar[r]
             & 0
           }
\]
The lower row is pure exact and the inclusion
$K^{\prime\prime}/K^{\prime} \subseteq F/K^{\prime}$ is pure, both by
construction.  Hence, if we tensor the diagram with an arbitrary
$R$-right-module $Q$ it follows from the snake lemma that $Q \otimes
K^{\prime\prime} \rightarrow Q \otimes F$ is injective.  So
$K^{\prime\prime}$ is a pure submodule of $F$ as desired.
\end{proof}

The following lemma is due to Angeleri-H\"{u}gel, Mantese, Tonolo, and
Trlifaj; cf.\ \cite[proof of lem.\ 9]{EklofTrlifaj}.

\begin{Lemma}
\label{lem:AHMTT}
If a class $\BF$ is closed under set indexed direct sums and pure
quotients, then it is also closed under colimits indexed by partially
ordered sets.
\end{Lemma}

\begin{proof}
Let $\{ F_i \}$ be a system in $\BF$ indexed by a partially ordered
set.  Then it is easy to see that the  canonical surjection
$\bigoplus F_i \rightarrow \dirlim F_i$ is a pure epimorphism, and the
lemma follows.
\end{proof}

\begin{Theorem}
\label{thm:covering}
If a class $\BF$ is closed under pure quotient modules, then the
following conditions are equivalent.
\begin{enumerate}

  \item  $\BF$ is closed under set indexed direct sums.

\smallskip

  \item  $\BF$ is precovering.

\smallskip

  \item  $\BF$ is covering.

\end{enumerate}
\end{Theorem}

\begin{proof}
Observe that $\BF$ is closed under direct summands, since the
projection onto a direct summand turns it into a pure quotient.
Moreover, $\BF$ is a cardinality class by Proposition
\ref{pro:Enochs}.

Proposition \ref{pro:precovering_preenveloping}(i) gives that (i) and
(ii) are equivalent.  By definition, (iii) implies (ii).

Finally, suppose that (ii) holds.  Then (i) also holds by the above,
and so Lemma \ref{lem:AHMTT} says that $\BF$ is closed under colimits
indexed by partially ordered sets, and in particular under well
ordered colimits.  But then $\BF$ is covering by \cite[thm.\
5.2.3]{EJ}, proving (iii).
\end{proof}

\begin{Remark}
\label{rmk:duals}
The dual of Proposition \ref{pro:Enochs} and the dual of Theorem
\ref{thm:covering} except the covering part were proved, up to
differences of terminology, by Rada and Saor\'{\i}n in \cite[prop.\ 
2.8 and cor.\ 3.5(c)]{RadaSaorin}.  Namely,
\begin{enumerate}
  
  \item If a class $\BG$ is closed under pure submodules, then it is a
        co-cardinality class (cf.\ Definition \ref{def:Enochs}(ii)).

\smallskip

  \item Let the class $\BG$ be closed under pure submodules.  Then
        $\BG$ is preenveloping if and only if it is closed under set
        indexed direct products.

\end{enumerate}

The covering part of Theorem \ref{thm:covering} cannot be dualized.
For example, if a ring is right-coherent, then the class of flat
left-modules is closed under set indexed products, and it is easy to
see that this class is also closed under pure submodules.  But flat
envelopes of left-modules do not necessarily exist, see \cite[thm.\ 
6.1]{Enochs}.
\end{Remark}

\section{Cotorsion pairs}
\label{sec:cotorsion}

This section shows that if $\BF$ is a class which contains the ground
ring and is closed under extensions, direct sums, pure submodules, and
pure quotients, then $\BF$ is the first half of a so-called perfect
cotorsion pair (see Definition \ref{def:cotorsion_pairs}).

Recall the following important notion from Enochs and L\'{o}pez-Ramos
\cite[def.\ 2.1]{ELR}.

\begin{Definition}
A class $\BF$ is called a {\em Kaplansky class} if there is a cardinal
number $\ff$ such that, when $F$ is in $\BF$ and $f$ is an element of
$F$, we have
\[
  f \in F^{\prime} \subseteq F
\]
for some submodule $F^{\prime}$ where $F^{\prime}$ and $F/F^{\prime}$
are in $\BF$ and where $|F^{\prime}| \leq \ff$.
\end{Definition}

\begin{Proposition}
\label{pro:Kaplansky}
If a class $\BF$ is closed under pure submodules and pure quotient
modules, then it is a Kaplansky class.
\end{Proposition}

\begin{proof}
Given $F$ in $\BF$ and $f$ in $F$, the submodule $Rf$ has $|Rf| \leq
|R|$.

By \cite[lem.\ 5.3.12]{EJ}, there is a cardinal number $\ff$ depending
only on $|R|$ such that we can enlarge $Rf$ to a pure submodule
$F^{\prime}$ of $F$ with $|F^{\prime}| \leq \ff$.  So
\[
  f \in F^{\prime} \subseteq F,
\]
and $F^{\prime}$ and $F/F^{\prime}$ are both in $\BF$ since $\BF$ is
closed under pure submodules and pure quotients.
\end{proof}

Let us recall the definition of a cotorsion pair.  This goes back to
Salce \cite{Salce}, and has gained popularity as a framework for
relative homological algebra; among our references we could mention
\cite{ATT}, \cite{BBE}, \cite{EklofTrlifaj}, \cite{EJ}, \cite{ELR},
and \cite{KS}.

\begin{Definition}
\label{def:cotorsion_pairs}
Let $\BF$ and $\BG$ be classes.  Then
\[
  \BF^{\perp} = \Ker \Ext^1(\BF,-)
              = \{\, N \in \Mod R \,|\, 
                    \Ext^1(F,N) = 0 \; \mbox{for} \; F \in \BF \,\}
\]
and
\[
  {}^{\perp}\BG = \Ker \Ext^1(-,\BG)
                = \{\, M \in \Mod R \,|\, 
                      \Ext^1(M,G) = 0 \; \mbox{for} \; G \in \BG \,\}.
\]
The pair $(\BF,\BG)$ is called a {\em cotorsion pair} if $\BF^{\perp}
= \BG$ and $\BF = {}^{\perp}\BG$.

The cotorsion pair is called {\em perfect} if $\BF$ is covering and
$\BG$ is enveloping. 
\end{Definition}

\begin{Theorem}
\label{thm:cotorsion}
If a class $\BF$ contains the ground ring $R$ and is closed under
extensions, set indexed direct sums, pure submodules, and pure
quotient modules, then $(\BF,\BF^{\perp})$ is a perfect cotorsion
pair.

In particular, $\BF$ is covering and $\BF^{\perp}$ is enveloping. 
\end{Theorem}

\begin{proof}
We shall use the powerful result \cite[thm.\ 2.9]{ELR} to prove this.  

To verify that the conditions of \cite[thm.\ 2.9]{ELR} are satisfied,
first note that $\BF$ is a Kaplansky class by Proposition
\ref{pro:Kaplansky}.

Since $\BF$ contains $R$ and is closed under set indexed direct sums,
it follows that $\BF$ contains all free modules.  But $\BF$ is closed
under pure quotients and so in particular under direct summands, and
so in fact, $\BF$ contains all projective modules.

Finally, since $\BF$ is closed under set indexed direct sums and pure
quotients, it is closed under all colimits indexed by partially
ordered sets by Lemma \ref{lem:AHMTT}.

This shows that the conditions of \cite[thm.\ 2.9]{ELR} are satisfied,
and the present theorem follows.
\end{proof}

\section{Applications}
\label{sec:applications}

This section gives a number of applications of the theory developed
above.

\begin{Remark}
\label{rmk:BC}
In the following results, note that the class
\[
  \Ker \Tor_1(\BB,-)
  = \{\, M \,|\, \Tor_1(B,M) = 0 \; \mbox{for each $B$ in $\BB$} \,\}
\]
can be obtained as ${}^{\perp}\BC$ by setting $\BC$ equal to the set
of all Pontryagin duals $\Hom_{\BZ}(B,\BQ/\BZ)$ for $B$ in $\BB$.
This holds by \cite[proof of cor.\ 11]{EklofTrlifaj} and also follows
from computation \eqref{equ:p} below.
\end{Remark}

\begin{Lemma}
\label{lem:ABC}
Let $\BC$ be a class of pure injective $R$-left-modules, $\BB$ a class
of $R$-right-modules, and $\BA$ a class of finitely presented
$R$-left-modules.
\begin{enumerate}

  \item  The class ${}^{\perp}\BC$ is closed under set indexed direct
         sums and pure quotients.

\smallskip

  \item  The class $\Ker \Tor_1(\BB,-)$ is closed under  set indexed
         direct sums, pure quotients, and pure submodules.

\smallskip
\noindent
         If $R$ is right-coherent and $\BB$ consists of finitely
         presented mo\-du\-les, then $\Ker \Tor_1(\BB,-)$ is closed
         under set indexed direct products.

  \item  The class $\BA^{\perp}$ is closed under
         set indexed direct products and pure submodules.

\smallskip
\noindent
         If $R$ is left-coherent, then $\BA^{\perp}$ is closed under
         set indexed direct sums and pure quotients.

\end{enumerate}
\end{Lemma}

\begin{proof}

\noindent
(i).  This is easy to prove, using the observation that for $C$ in
$\BC$, the functor $\Hom(-,C)$ sends pure exact sequences to exact
sequences. 

\noindent
(ii).  Since $\Ker \Tor_1(\BB,-)$ has the form ${}^{\perp}\BC$ for a
suitable set $\BC$ by Remark \ref{rmk:BC}, the statements about set
indexed direct sums and pure quotients follow from part (i).

If
\[
  0 \rightarrow X^{\prime}
    \rightarrow X
    \rightarrow X^{\prime\prime}
    \rightarrow 0
\]
is a pure exact sequence, then by \cite[thm.\ 6.4]{JL}, the
Pontryagin duality functor $(-)^{\vee} = \Hom_{\BZ}(-,\BQ/\BZ)$
gives a split exact sequence
\[
  0 \rightarrow (X^{\prime\prime})^{\vee}
    \rightarrow X^{\vee}
    \rightarrow (X^{\prime})^{\vee}
    \rightarrow 0,
\]
so if $B$ is in $\BB$ then there is a split exact sequence
\begin{equation}
\label{equ:c}
  0 \rightarrow \Ext_{R^{\opp}}^1(B,(X^{\prime\prime})^{\vee})
    \rightarrow \Ext_{R^{\opp}}^1(B,X^{\vee})
    \rightarrow \Ext_{R^{\opp}}^1(B,(X^{\prime})^{\vee})
    \rightarrow 0.
\end{equation}
Moreover, a standard computation shows
\begin{align}
\nonumber
  \Ext_{R^{\opp}}^1(B,(-)^{\vee}) &
  = \Ext_{R^{\opp}}^1(B,\Hom_{\BZ}(-,\BQ/\BZ)) \\
\nonumber
  & \simeq \Hom_{\BZ}(\Tor^R_1(B,-),\BQ/\BZ) \\
\label{equ:p}
  & = \Tor^R_1(B,-)^{\vee}.
\end{align}
Now let $X$ be in $\Ker \Tor_1(\BB,-)$ so $\Tor^R_1(B,X) = 0$.  The
last computation implies $\Ext_{R^{\opp}}^1(B,X^{\vee}) = 0$.  The
sequence \eqref{equ:c} shows
\[
  \Ext_{R^{\opp}}^1(B,(X^{\prime})^{\vee}) = 0,
\]
and then the last computation again implies that
$\Tor^R_1(B,X^{\prime}) = 0$.  So $X^{\prime}$ is in $\Ker
\Tor_1(\BB,-)$.

Finally, if $R$ is right-coherent and $\BB$ consists of finitely
presented mo\-du\-les, then each $B$ in $\BB$ has a projective
resolution consisting of finitely generated modules, so $\Ker
\Tor_1(\BB,-)$ is closed under set indexed direct products because
these are preserved by the functor $\Tor_1(B,-)$.

(iii).  It is clear that $\BA^{\perp}$ is closed under set indexed
direct pro\-ducts because these are preserved by the functor
$\Ext^1(A,-)$.

If
\begin{equation}
\label{equ:y}
  0 \rightarrow Y^{\prime}
    \rightarrow Y
    \rightarrow Y^{\prime\prime}
    \rightarrow 0
\end{equation}
is a pure exact sequence and $A$ is in $\BA$, then there is an exact
sequence 
\[
  \Hom(A,Y)
  \rightarrow \Hom(A,Y^{\prime\prime})
  \rightarrow \Ext^1(A,Y^{\prime})
  \rightarrow \Ext^1(A,Y).
\]
The first arrow is surjective because $A$ is finitely presented, so
the second arrow is zero.  If $Y$ is in $\BA^{\perp}$ then
$\Ext^1(A,Y) = 0$, but then the sequence shows $\Ext^1(A,Y^{\prime}) =
0$ whence $Y^{\prime}$ is in $\BA^{\perp}$.

Now suppose that $R$ is left-coherent.  By \cite[thm.\ 6.4]{JL} again,
the Pontryagin duality functor $(-)^{\vee} = \Hom_{\BZ}(-,\BQ/\BZ)$
sends the pure exact sequence \eqref{equ:y} to a split exact sequence
\[
  0 \rightarrow (Y^{\prime\prime})^{\vee}
    \rightarrow Y^{\vee}
    \rightarrow (Y^{\prime})^{\vee}
    \rightarrow 0,
\]
so if $A$ is in $\BA$ then there is a split exact sequence
\begin{equation}
\label{equ:a}
  0 \rightarrow \Tor_1((Y^{\prime\prime})^{\vee},A)
    \rightarrow \Tor_1(Y^{\vee},A)
    \rightarrow \Tor_1((Y^{\prime})^{\vee},A)
    \rightarrow 0.
\end{equation}
However, $A$ has a projective resolution consisting of finitely
generated modules, so a standard computation shows
\begin{align*}
  \Tor^R_1((-)^{\vee},A) & = \Tor^R_1(\Hom_{\BZ}(-,\BQ/\BZ),A) \\
  & \simeq \Hom_{\BZ}(\Ext_R^1(A,-),\BQ/\BZ) \\
  & = \Ext_R^1(A,-)^{\vee}.
\end{align*}
Now let $Y$ be in $\BA^{\perp}$ so $\Ext^1(A,Y) = 0$.  The last
computation implies $\Tor_1(Y^{\vee},A) = 0$.  The sequence
\eqref{equ:a} shows that $\Tor_1((Y^{\prime\prime})^{\vee},A) = 0$,
and then the last computation again implies that
$\Ext^1(A,Y^{\prime\prime}) = 0$.  So $Y^{\prime\prime}$ is in
$\BA^{\perp}$.

Finally, since each $A$ in $\BA$ has a projective resolution
consisting of finitely generated modules, $\BA^{\perp}$ is closed
under set indexed direct sums because these are preserved by the
functor $\Ext^1(A,-)$.
\end{proof}

Some parts of the following theorem were already known; for instance,
it was proved by Eklof and Trlifaj in \cite[cor.\ 10 and thm.\ 
12(i)]{EklofTrlifaj} that $\Ker \Tor_1(\BB,-)$ and ${}^{\perp}\BC$ are
both covering, and when the ground ring is left-coherent, El Bashir's
result \cite[thm.\ 3.3]{ElBashir} implies that $\BA^{\perp}$ is also
covering if one is willing to assume Vopenka's Principle on high
cardinal numbers.  However, it is new that we are able to work with
the intersections of such classes.

\begin{Theorem}
\label{thm:ABC}
Let $\BC$ be a class of pure injective $R$-left-modules, $\BB$ a class
of $R$-right-modules, and $\BA$ a class of finitely presented
$R$-left-modules.
\begin{enumerate}
  
  \item  The classes
\[
  {}^{\perp}\BC \;\;\;\;\mbox{and}\;\;\;\;  \Ker \Tor_1(\BB,-)
\]
         are covering.

\smallskip
\noindent
         If $R$ is left-coherent, then the classes
\[
  \BA^{\perp} \textstyle{\bigcap} {}^{\perp}\BC \;\;\;\;\mbox{and}\;\;\;\;
  \BA^{\perp} \bigcap \Ker \Tor_1(\BB,-)
\]
         are covering.

\smallskip

  \item  The class
\[
  \Ker \Tor_1(\BB,-)
\]
         is the first half of  a perfect cotorsion pair.

\smallskip
\noindent
         If $R$ is left-coherent and is contained in $\BA^{\perp}$,
         then the class 
\[
  \BA^{\perp} \textstyle{\bigcap} \Ker \Tor_1(\BB,-)
\]
         is the first half of a perfect cotorsion pair.

\smallskip

  \item  The class
\[
  \BA^{\perp}
\]
         is preenveloping.

\smallskip
\noindent
         If $R$ is right-coherent and $\BB$ consists of finitely
         presented mo\-du\-les, then the class
\[
  \BA^{\perp} \textstyle{\bigcap} \Ker \Tor_1(\BB,-)
\]
         is preenveloping.
\end{enumerate}
\end{Theorem}

\begin{proof}
(i).  It is enough to prove the statements involving ${}^{\perp}\BC$, 
since $\Ker \Tor_1(\BB,-)$ has the form ${}^{\perp}\BC$ 
by Remark \ref{rmk:BC}.

The class ${}^{\perp}\BC$ is closed under set indexed direct
sums and pure quotients by Lemma \ref{lem:ABC}(i).  If $R$ is
left-coherent, then $\BA^{\perp}$ has the same properties by Lemma
\ref{lem:ABC}(iii).  Now use Theorem \ref{thm:covering}.

\noindent
(ii).  The class $\Ker \Tor_1(\BB,-)$ clearly contains $R$ and
is closed under extensions, and by Lemma \ref{lem:ABC}(ii) it is also
closed under set indexed direct sums, pure quotients, and pure
submodules.  If $R$ is left-coherent and is contained in
$\BA^{\perp}$, then $\BA^{\perp}$ has the same properties by Lemma
\ref{lem:ABC}(iii), and so $\BA^{\perp} \bigcap \Ker \Tor_1(\BB,-)$
also has the same properties.  Now use Theorem \ref{thm:cotorsion}.

\noindent
(iii).  The class $\BA^{\perp}$ is closed under set indexed direct
products and pure submodules by Lemma \ref{lem:ABC}(iii).  If $R$ is
right-coherent and $\BB$ consists of finitely presented modules, then
$\Ker \Tor_1(\BB,-)$ has the same properties by Lemma
\ref{lem:ABC}(ii), and so $\BA^{\perp} \bigcap \Ker \Tor_1(\BB,-)$ has
the same properties.  Now use Remark \ref{rmk:duals}(ii).
\end{proof}

\begin{Example}
\label{exa:ABC}
\begin{enumerate}

  \item  Let $m$ be an non-negative integer and consider the class
\[
  \BF_{\leq m} = \{\, F \,|\, \mbox{$F$ is an $R$-left-module with
                                    flat dimension $\leq m$} \,\}.
\]
Then $(\BF_{\leq m},(\BF_{\leq m})^{\perp})$ is a perfect cotorsion
pair.  In particular, $\BF_{\leq m}$ is covering.

\smallskip
\noindent
Moreover, if $R$ is right-coherent, then $\BF_{\leq m}$ is also
pre\-en\-ve\-lo\-ping.

\smallskip
\noindent
This follows from Theorem \ref{thm:ABC}, (ii) and (iii), by setting
$\BB$ equal to the $m$'th syzygies in projective resolutions of
finitely presented modules; cf.\ \cite[thm.\ A.8]{JL}.

\smallskip

  \item   Suppose that $R$ is left-noetherian and let $n$ be an
non-negative integer.  Then the class 
\[
  \BI_{\leq n} = \{\, I \,|\, \mbox{$I$ is an $R$-left-module with
                                    injective dimension $\leq n$} \,\}
\]
is covering and preenveloping.

\smallskip
\noindent
This follows from Theorem \ref{thm:ABC}, (i) and (iii), by setting
$\BA$ equal to the $n$'th syzygies in projective resolutions of
finitely ge\-ne\-ra\-ted modules; cf.\ \cite[thm.\ A.6]{JL}.

\smallskip

  \item  Suppose that $R$ is left-noetherian and let $m$ and $n$ be
    non-negative integers.  Then the class
\[
  \BF_{\leq m} \cap \BI_{\leq n} = \biggl\{\, X \,\bigg|\, 
  \begin{array}{l}
    \mbox{$X$ is an $R$-left module with flat dimen-} \\
    \mbox{sion $\leq m$ and injective dimension $\leq n$} 
  \end{array}
                 \,\biggl\}
\]
is covering.  

\smallskip
\noindent
Moreover, if $R$ is right-coherent, then $\BF_{\leq m} \cap \BI_{\leq
n}$ is also preenveloping. 

\smallskip
\noindent
This follows from Theorem \ref{thm:ABC}, (i) and (iii), using the same
$\BA$ and $\BB$ as above.

\smallskip

  \item  The class of fp-injective $R$-left-modules,
\[
  \BJ = \{\, J \,|\, \Ext^1(A,J) = 0 \;\mbox{for $A$ a finitely
                     presented $R$-left-module} \,\},
\]
         is preenveloping. 

\smallskip
\noindent
         Moreover, if $R$ is left-coherent, then $\BJ$ is also
         covering.

\smallskip
\noindent
         This follows from Theorem \ref{thm:ABC}, (i) and (iii), by
         setting $\BA$ equal to all the finitely presented modules.

\smallskip
\noindent
         These results on $\BJ$ were known from \cite{Adams} (see also
         \cite[prop.\ 6.2.4]{EJ}) and \cite{Pinzon} with different
         pro\-ofs. 

\end{enumerate}
\end{Example}

Finally, we use our methods to give new proofs of some known results
about the torsion free modules in a torsion pair.

\begin{Definition}
\label{def:torsion_pairs}
Recall from \cite{Dickson} that a pair of classes $(\BT,\BF)$ is
called a {\em torsion pair} if $\BT \cap \BF$ contains only modules
isomorphic to $0$, the class $\BT$ is closed under quotient modules,
the class $\BF$ is closed under submodules, and each module $M$
permits a short exact sequence
\[
  0 \rightarrow T \rightarrow M \rightarrow F \rightarrow 0
\]
with $T$ in $\BT$ and $F$ in $\BF$.

The torsion pair is called {\em hereditary} if $\BT$ is also closed
under sub\-mo\-du\-les, see \cite[p.\ 441]{Teply2}.

The torsion pair is said to be of {\em finite type} if each left-ideal
$\fa$ for which $R/\fa$ is in $\BT$ contains a finitely generated
left-ideal $\fb$ for which $R/\fb$ is in $\BT$, see \cite[p.\
649]{BT}.  Note that if $R$ is left-noetherian, then the torsion pair
is automatically of finite type.
\end{Definition}

\begin{Lemma}
\label{lem:hereditarytorsion}
Let $(\BT,\BF)$ be a hereditary torsion pair.
\begin{enumerate}

  \item  $F$ is in $\BF$ if and only if 
\[
  \Hom(R/\fa,F) = 0
\]
         for each ideal $\fa$ such that $R/\fa$ is in $\BT$.

\smallskip

  \item  If $(\BT,\BF)$ is of finite type, then $F$ is in
         $\BF$ if and only if
\[
  \Hom(R/\fb,F) = 0
\]
         for each finitely generated left-ideal $\fb$ such that
         $R/\fb$ is in $\BT$.
\end{enumerate}
\end{Lemma}

\begin{proof}
The module $F$ is in $\BF$ if and only if $\Hom(T,F) = 0$ for each
$T$ in $\BT$, see \cite{Dickson}.  It is a small computation to see
that this implies the lemma's statements.
\end{proof}

\begin{Lemma}
\label{lem:properties_torsion_theory}
Let $(\BT,\BF)$ be a torsion pair.  Then
\begin{enumerate}

  \item  $\BF$ is closed under extensions.

\smallskip

  \item  If $(\BT,\BF)$ is hereditary, then $\BF$ is closed under set
         indexed direct sums. 

\smallskip

  \item  If $(\BT,\BF)$ is hereditary and of finite type, then $\BF$ is
         closed under pure submodules and pure quotient modules.

\end{enumerate}
\end{Lemma}

\begin{proof}
(i).  This holds because $F$ is in $\BF$ if and only if
\[
  \Hom(T,F) = 0
\]
for each $T$ in $\BT$.

(ii).  Let $\{ F_i \}$ be a set indexed system in $\BF$.  Let $\fa$ be
a left-ideal in $R$ with $R/\fa$ in $\BT$.  Since $R/\fa$ is finitely
generated, we get the following $\cong$,
\[
  \Hom(R/\fa,\bigoplus F_i) \cong \bigoplus \Hom(R/\fa,F_i) = 0,
\]
where $=$ is because each $F_i$ is in $\BF$.  By Lemma
\ref{lem:hereditarytorsion}(i) this shows that $\bigoplus F_i$ is in
$\BF$.

\smallskip

(iii).  Let $F$ be in $\BF$ and let 
\[
  0 \rightarrow F^{\prime} \rightarrow F \rightarrow F^{\prime\prime}
  \rightarrow 0
\]
be a pure exact sequence.  As $\BF$ is closed under submodules,
$F^{\prime}$ is in $\BF$ as desired.

Let $\fb$ be a finitely generated left-ideal in $R$ with $R/\fb$ in
$\BT$, and let $R/\fb \rightarrow F^{\prime\prime}$ be a homomorphism.
Since $R/\fb$ is finitely presented, $R/\fb \rightarrow
F^{\prime\prime}$ factors through the pure epimorphism $F \rightarrow
F^{\prime\prime}$.  But $F$ is in $\BF$ so each homomorphism $R/\fb
\rightarrow F$ is zero, and it follows that $R/\fb \rightarrow
F^{\prime\prime}$ is zero.  Hence $F^{\prime\prime}$ is in $\BF$ by
Lemma \ref{lem:hereditarytorsion}(ii), as desired.
\end{proof}

The following result was first proved by Bican and Torrecillas in
\cite[cor.\ 4.1]{BT2}.

\begin{Theorem}
\label{thm:torsion_free_covering}
Let $(\BT,\BF)$ be a hereditary torsion pair of finite type.

Then $\BF$ is covering.
\end{Theorem}

\begin{proof} 
Lemma \ref{lem:properties_torsion_theory}(ii) says that $\BF$ is
closed under set indexed direct sums, and Lemma
\ref{lem:properties_torsion_theory}(iii) says that $\BF$ is closed
under pure quotients, so $\BF$ is covering by Theorem
\ref{thm:covering}.
\end{proof}

The following result was first proved by Angeleri-H\"{u}gel, Tonolo,
and Trlifaj in \cite[exa.\ 2.7]{ATT}.

\begin{Theorem}
\label{thm:torsion_and_cotorsion}
Let $(\BT,\BF)$ be a hereditary torsion pair of finite type where the
ground ring $R$ is in $\BF$.

Then $(\BF,\BF^{\perp})$ is a perfect cotorsion pair.

In particular, $\BF$ is covering and $\BF^{\perp}$ is enveloping.
\end{Theorem}

\begin{proof}
Lemma \ref{lem:properties_torsion_theory} says that $\BF$ is
closed under extensions, set indexed direct sums, pure submodules,
and pure quotients.  As $R$ is in $\BF$, it follows that
$(\BF,\BF^{\perp})$ is a perfect cotorsion pair by Theorem
\ref{thm:cotorsion}.
\end{proof}

\medskip
{\em Acknowledgement. }
We thank Lidia Angeleri-H\"{u}gel, Ladislav Bican, Robert El Bashir,
Edgar E.\ Enochs, Juan Antonio L\'{o}pez-Ramos, and Jan Trlifaj warmly
for numerous comments to the preliminary versions of this paper.
Their expert advice has led to substantial improvements.

We thank Katherine R.\ Pinzon for com\-mu\-ni\-ca\-ting
\cite{Pinzon}.

The second author was supported by a grant from the Royal Society.

\end{document}